\documentclass[12pt]{amsart}
\input amssym.def

\newcommand{\lr}{\mbox{$\longrightarrow$}}

\newcommand{\be}{\begin{equation}}
\newcommand{\ee}{\end{equation}}
\newtheorem{guess}{Theorem}[section]
\newcommand{\bth}{\begin{guess}$\!\!\!$}
\newcommand{\eeth}{\end{guess}}
\newtheorem{propo}[guess]{Proposition}
\newcommand{\bpropo}{\begin{propo}$\!\!$~}
\newcommand{\epropo}{\end{propo}}
\newtheorem{lema}[guess]{Lemma}
\newcommand{\blem}{\begin{lema}$\!\!\!$}
\newcommand{\elem}{\end{lema}}
\newtheorem{defe}[guess]{Definition}
\newcommand{\bdefe}{\begin{defe}$\!\!\!$}
\newcommand{\edefe}{\end{defe}}
\newtheorem{coro}[guess]{Corollary}
\newcommand{\bcor}{\begin{coro}$\!\!\!$}
\newcommand{\ecor}{\end{coro}}
\newtheorem{rema}[guess]{Remark}
\newcommand{\brem}{\begin{rema}$\!\!\!$~\rm}
\newcommand{\erem}{\end{rema}}
\newtheorem{exam}[guess]{Example}
\newcommand{\beg}{\begin{exam}$\!\!\!$~\rm}
\newcommand{\eeg}{\end{exam}}

\newcommand{\bz}{{\Bbb Z}}
\newcommand{\bc}{{\Bbb C}}

\newcommand{\bq}{{\Bbb  Q}}

\newcommand{\br}{{\Bbb R}}
\newcommand{\bp}{{\Bbb P}}

\renewcommand{\phi}{{\varphi}}

\newcommand{\bn}{{\Bbb N}}
\begin{document}
\title{Degrees of maps between Grassmann manifolds}
\author[P. Sankaran]{Parameswaran Sankaran}
\address{The Institute of Mathematical Sciences, Chennai 600113.}
\email{sankaran@imsc.res.in}

\author[S. Sarkar]{Swagata Sarkar}
\address{The Institute of Mathematical Sciences, Chennai 
600113.}
\email{swagata@imsc.res.in}
\thispagestyle{empty}
\maketitle

\noindent
{\bf Abstract:}{\it Let   $f\colon {\Bbb G}_{n,k}\lr {\Bbb G}_{m,l}$ be any continuous map between two {\it distinct} complex 
(resp. quaternionic) Grassmann manifolds of the same dimension.  We show that the degree of $f$ is zero provided $n,m$ are sufficiently large and $l\geq 2$.  If the degree of $f$ is $\pm 1$, we show that $(m,l)=(n,k)$ and  $f$ is a homotopy equivalence.  Also, we prove that the image under $f^*$ of elements of a set of algebra generators of $H^*({\Bbb G}_{m,l};\bq)$ is determined  upto a sign, $\pm$, if the degree of $f$ is non-zero.}

\section{Introduction}
The purpose of this paper is to study degrees of 
maps between two distinct complex (resp.  quaternionic)  Grassmann manifolds. It 
can be viewed as a continuation of  the paper \cite{rs} 
where the case of oriented (real) Grassmann manifolds was 
settled completely.  The same problem in the case of complex and 
quaternionic Grassmann manifolds was initiated and settled in \cite{rs} 
in half the cases.   
The problem can be formulated 
purely algebraically in terms of algebra homomorphism between 
the cohomology algebras of the complex Grassmann manifolds concerned.  These algebras have additional structures, arising from Poincar\'e duality and the hard Lefschetz theorem. Our results are obtained by exploiting these properties. 
In view of the fact that the integral cohomology ring of 
a quaternionic Grassmann manifold  
is isomorphic to that of the corresponding complex Grassmann 
manifold via a degree doubling isomorphism, 
and since our proofs involve mostly analyzing the algebra- homomorphisms between the cohomology algebras 
of the Grassmann manifolds, 
we will only need to consider the case of complex 
Grassmann manifolds.  (In 
the course of our proof of Theorem \ref{degreeone},   
simply-connectedness of the complex Grassmann 
manifold will be used; the same property 
also holds for the quaternionic Grassmann manifolds.)
For this reason, we need only to consider the case of complex Grassmann manifolds. 
 
Let 
${\Bbb F}$ denote the field $\bc$ of complex numbers 
or the skew-field ${\Bbb H}$ of quaternions. We denote by ${\Bbb F}{\Bbb G}_{n,k}$ the ${\Bbb F}$-Grassmann manifold of $k$-dimensional 
left ${\Bbb F}$-vector subspaces of ${\Bbb F}^n$. 
Let $d:=\dim_{\Bbb R}{\Bbb F}$. 
 Since we will mostly deal with complex Grassmann manifolds, we shall write ${\Bbb G}_{n,k}$ 
instead of ${\Bbb C}{\Bbb G}_{n,k}$;  the phrase 
`Grassmann manifold',  without further qualification, will always mean a complex Grassmann manifold. 

Using the usual `hermitian' metric on ${\Bbb F}^n$, one 
obtains a diffeomorphism $\perp\colon {\Bbb F}{\Bbb G}_{n,k}\cong {\Bbb F}{\Bbb G}_{n,n-k}$. For this reason, it suffices to consider 
only those ${\Bbb F}$-Grassmann manifolds ${\Bbb F}{\Bbb G}_{n,k} $ with $1\leq k\leq 
[n/2]$.  Let $1\leq l\leq [m/2]$ be another ${\Bbb F}$-Grassmann manifold 
having the same dimension as ${\Bbb F}{\Bbb G}_{n,k}$ so that 
$\dim_{{\Bbb F}}{\Bbb F}{\Bbb G}_{n,k}=k(n-k)=l(m-l)=:N$.  

Complex Grassmann manifolds admit a 
natural orientation arising from the 
fact they have a natural complex structure. Although 
the quaternionic Grassmann manifolds do not admit 
even almost complex structures (cf. \cite{markl}), 
they are simply connected and hence orientable.    

Let $f\colon {\Bbb F}{\Bbb G}_{n,k}\lr {\Bbb F}{\Bbb G}_{m,l}$ be any continuous map.  
It was observed in \cite{rs} that when $1\leq k<l\leq 
[m/2]$, the degree of $f$ is zero.  When $l=1$, one has $N=m-1$ and  
${\Bbb F}{\Bbb G}_{m,l}$ is just the ${\Bbb F}$- projective space ${\Bbb F}\bp^{N}$. The set of  
homotopy classes of maps $f\colon {\Bbb F}{\Bbb G}_{n,k}\lr {\Bbb F}\bp^{N}$ 
are in bijection with homomorphisms of abelian groups 
$\bz\cong 
H^{d}({\Bbb F}\bp^{N};\bz)\lr H^{d}({\Bbb F}{\Bbb G}_{n,k};\bz)\cong \bz$  where $d=\dim_\br{\Bbb F}$, via the induced homomorphism.  Furthermore 
the degree of $f$ is determined by $f^*\colon H^d({\Bbb F}\bp^{N};\bz)\lr H^d({\Bbb F}{\Bbb G}_{n,k};\bz)$. 
(See \cite{rs} for details.)

We now state the main results of this paper.
 
\bth \label{degdet} Let ${\Bbb F}=\bc$ or $ {\Bbb H}$ 
and let $d=\dim_\br {\Bbb F}$.
Let $f\colon {\Bbb F}{\Bbb G}_{n,k}\lr {\Bbb F}{\Bbb G}_{m,l}$ be any continuous map between 
two ${\Bbb F}$-Grassmann manifolds of the same dimension. Then, there exist   
algebra generators 
$u_i\in H^{di}({\Bbb F}{\Bbb G}_{m,l};\bq), 1\leq i\leq l,$
such that the image 
$f^*(u_i)\in H^{di}({\Bbb F}{\Bbb G}_{n,k};\bq)$, $1\leq i\leq l,$ is 
determined upto  a sign $\pm$, provided degree of $f$ is non-zero. 
\eeth

\bth \label{degreezero} Let ${\Bbb F}=\bc$ or ${\Bbb H}$. 
Fix integers $2\leq l<k$. Let $m, n\geq 2k$ be positive integers 
such that $k(n-k)=l(m-l)$ and $f\colon {\Bbb F}{\Bbb G}_{n,k}\lr {\Bbb F}{\Bbb G}_{m,l}$ any continuous map. Then, degree of $f$ is zero if 
$(l^2-1)(k^2-1)((m-l)^2-1)((n-k)^2-1)$ is not a perfect 
square.  In particular, degree of  $f$ is zero for $n$ 
sufficiently large.  
\eeth 

\bth \label{degreeone} Let ${\Bbb F}=\bc$ or ${\Bbb H}$. 
Suppose that $k(n-k)=l(m-l),$ and $1\leq l\leq [m/2], 1\leq k\leq [n/2]$. 
If $f\colon {\Bbb F}{\Bbb G}_{n,k}\lr  {\Bbb F}{\Bbb G}_{m,l}$ is a map of degree $\pm 1$, then 
$(m,l)=(n,k)$ and $f$ is a homotopy equivalence.
\eeth

Our proofs make use of the 
notion of degrees of Schubert varieties, extended to cohomology classes. Theorem \ref{degreeone}, which is 
an analogue in the topological realm of a result of K. H. Paranajape and V. Srinivas \cite{parsri}, is proved using Whitehead's theorem.  
Proof of Theorem \ref{degdet} uses some properties of 
the cohomology of the  complex Grassmann manifolds arising from Hodge theory. (See \ref{mainpropo}.)  Theorem \ref{degreezero} is proved by reducing it to a diophantine problem and 
appealing  to  Siegel's Theorem on solutions of certain polynomial equation of the form $y^2=F(x)$. In 
our situation, $F(x)$ will be of degree $4$ over $\bq$ 
having distinct zeros.  

We now highlight  the following conjecture  
made in \cite{rs}. Theorem \ref{degreezero} provides the strongest evidence in support of the conjecture. 

\noindent 
{\bf Conjecture:} Let ${\Bbb F}=\bc$ or ${\Bbb H}$ and 
let $2\leq l<k\leq n/2<m/2$ where $k,l,m,n\in \bn.$  
Assume that $k(n-k)=l(m-l)$. Let $f\colon {\Bbb F}{\Bbb G}_{n,k}\lr {\Bbb F}{\Bbb G}_{m,l}$ 
be {\it any} continuous map. The degree of $f$ is zero. 

The paper is organized as follows. In \S 2 we recall 
basic and well-known facts concerning the cohomology algebra of the complex Grassmann manifolds. We shall 
consider continuous maps from a cohomologically 
K\"ahler manifold and establish some important 
properties in \S 3. They will be used in the course of our proofs. We prove the above theorems in \S 4. 

\section{Cohomology of Grassmann manifolds}
There are at least two well-known descriptions of the cohomology ring of a complex Grassmann manifold 
${\Bbb G}_{n,k}$.  We recall both of them. 

Let $\gamma_{n,k}$ be the `tautological' bundle over ${\Bbb G}_{n,k}$ whose fibre over a point $V\in {\Bbb G}_{n,k}$ is the 
$k$-dimensional complex vector space $V$. Evidently 
$\gamma_{n,k}$ is rank $k$-subbundle of the rank $n$ trivial bundle
$\mathcal{E}^n$ with projection $pr_1\colon {\Bbb G}_{n,k}\times \bc^n\lr {\Bbb G}_{n,k}$.
The quotient bundle $\mathcal{E}^n/\gamma_{n,k}$ 
is isomorphic to the orthogonal complement $\gamma_{n,k}^\perp$ 
in $\mathcal{E}^n$ 
(with respect to a hermitian metric on $\bc^n$) of the bundle $\gamma_{n,k}$. Let $c_i(\gamma_{n,k})\in H^{2i}({\Bbb G}_{n,k};\bz)$, 
be the $i$-th Chern class of $\gamma_{n,k}, 1\leq i\leq k$. Denoting the total Chern class of a vector bundle 
$\eta$  by $c(\eta) $ we see that $c(\gamma_{n,k}).c(\gamma_{n,k}^\perp)=1$. 

Let $c_1,\cdots, c_k$ denote the elementary symmetric 
polynomials in $k$ indeterminates $x_1,\cdots,x_k$.
Define  $h_j=h_j(c_1,\cdots, c_k)$ by the identity 
$$\prod_{1\leq i\leq k}(1+x_it)^{-1}=\sum_{j\geq 0}h_jt^j.$$    

Thus $c_j(\gamma^\perp_{n,k})=h_j(c_1(\gamma_{n,k}),c_2(\gamma_{n,k}),\cdots, c_k(\gamma_{n,k})), ~
1\leq j\leq n-k.$ (See \cite{ms}.)

Consider the ring $\bz[c_1,\cdots, c_k]/\mathcal{I}_{n,k}$ where degree of  $c_i=2i$, and $\mathcal{I}_{n,k}$ 
is the ideal $\langle h_j\mid j>n-k\rangle$. It can be shown that the elements $h_j, ~n-k+1\leq j\leq n,$  generate $\mathcal{I}_{n,k}.$ 
The homomorphism of graded rings   
$\bz[c_1,\cdots,c_k]\lr H^*({\Bbb G}_{n,k};\bz)$ defined by $c_i\mapsto c_i(\gamma_{n,k})$ is surjective and 
has kernel $\mathcal{I}_{n,k}$ and hence we have 
an isomorphism $H^*({\Bbb G}_{n,k};\bz)\cong
\bz[c_1,\cdots, c_k]/\mathcal{I}_{n,k}$.  
Henceforth we shall write $c_i$ to mean $c_i(\gamma_{n,k})\in H^*({\Bbb G}_{n,k};\bz)$. We shall denote 
by  $\bar{c}_j$ the element $c_j(\gamma^\perp_{n,k})
=h_j\in H^{2j}({\Bbb G}_{n,k};\bz)$.

As an abelian group, 
$H^*({\Bbb G}_{n,k};\bz)$ is free  
of rank ${n}\choose{k}$.  A $\bq$-basis for 
$H^{2r}({\Bbb G}_{n,k};\bq)$ is the set $\mathcal{C}_r$ of all 
monomials $c_1^{j_1} 
\cdots c_k^{j_k}$ where $j_i\leq n-k~\forall i, \sum_{1\leq i\leq k}ij_i=r$.  
In particular,  $c_k^{n-k}$ generates $H^{2N}({\Bbb G}_{n,k};\bq)\cong \bq$. If  ${\bf j}$ denotes the sequence $j_1,\cdots, j_k$, we shall denote by 
$c^{{\bf j}}$  the monomial $c_i^{j_1}\cdots c_k^{j_k}$. 
If $k\leq n/2$,   the 
set $\bar{\mathcal {C}}_r:=\{\bar{c}^{{\bf j}}\mid c^{{\bf j}}\in 
\mathcal{C}_r\}$ is also a basis for 
$H^{2r}({\Bbb G}_{n,k};\bq)$ 
where $\bar{c}^{{\bf j}}:=\bar{c}_1^{j_1}\cdots \bar{c}_k^{j_k}$. 

\noindent
{\bf Schubert calculus}\\
Another, more classical description of the cohomology 
ring of the Grassmann manifold ${\Bbb G}_{n,k}$ is 
via the Schubert calculus. Recall that ${\Bbb G}_{n,k}=\textrm{SL}(n,\bc)/P_k$ for the parabolic subgroup $P_k\subset \textrm{SL}(n,\bc)$ which stabilizes $\bc^k\subset \bc^n$ 
spanned by $e_1,\cdots, e_k$; here $e_i, 1\leq i\leq n,$ 
are the standard basis elements 
of $\bc^n$. Denote by $B\subset \textrm{SL}(n,\bc)$ the Borel subgroup of $\textrm{SL}(n,\bc)$  
which preserves the flag $\bc^1\subset \cdots\subset 
\bc^n$ and by $T\subset B$ the maximal torus which 
preserves the coordinate axes $\bc e_j, 1\leq j\leq n$. 
Let $I(n,k)$ denote the set of all $k$ element subsets of $\{1,2, \cdots, n\}$; we regard elements of $I(n,k)$ as increasing sequences of positive integers 
${\bf i}:= i_1<\cdots <i_k$ where $i_k\leq n$.  One has a partial order 
on $I(n,k)$ where, by definition, ${\bf i}\leq {\bf j}$ if $i_p\leq j_p$ 
for all $p, 1\leq p\leq k$.   
Let ${\bf i}\in I(n,k)$ and let $E_{{\bf i}}\in {\Bbb G}_{n,k}$ denote the vector 
subspace  of $\bc^n$ spanned by $\{e_j\mid j\in {\bf i}\}$.  The fixed points for the action of $T\subset \textrm{SL}(n)$ 
on ${\Bbb G}_{n,k}$ are precisely the $E_{{\bf i}}$, ${\bf i}
\in I(n,k)$.   

Schubert varieties in ${\Bbb G}_{n,k}$  are in bijection with the set $I(n,k)$. 
The $B$-orbit  of the $T$-fixed point $E_{{\bf i}}$ is the Schubert 
cell corresponding to ${\bf i}$ and is  
isomorphic to the affine space of (complex) dimension 
$\sum_j (i_j-j)=:|{\bf i}|$; its closure, denoted $\Omega_{{\bf i}},$ is the Schubert variety corresponding to ${\bf i}\in I(n,k)$. It is the union 
of all Schubert cells corresponding to those ${\bf j}\in I(n,k)$ such that  
${\bf j}\leq {\bf i}$.   
Schubert cells yield a cell decomposition of ${\Bbb G}_{n,k}$. Since the cells have even (real) dimension,  
the class of Schubert varieties form a $\bz$-basis for the 
integral homology of ${\Bbb G}_{n,k}$. Denote by  $[\Omega_{{\bf i}}]\in H^{2(N-|{\bf i}|)}({\Bbb G}_{n,k};\bz)$ the fundamental {\it dual} cohomology  
class determined by $\Omega_{{\bf i}}$.  (Thus $[{\Bbb G}_{n,k}]\in H^0({\Bbb G}_{n,k};\bz)$ is the identity element of the cohomology ring.)  
We shall denote the fundamental homology 
class of ${\Bbb G}_{n,k}$ by $\mu_{n,k}\in H^{2N}({\Bbb G}_{n,k};\bz)$.  

Schubert varieties corresponding to $
(n-k+1-i,n-k+2,\cdots, n)\in I(n,k), 0\leq i\leq n-k,$ are called {\it special} and 
will be denoted $\Omega_i$.  
More generally, if $\nu=\nu_1\geq \cdots
\geq \nu_k\geq 0$ is a partition of an integer $r$, $0\leq r\leq N$, 
with $\nu_1\leq n-k$, we obtain an element $(n-k+1-\nu_1,n-k+2-\nu_2,\cdots,n-\nu_k) \in I(n,k)$ with $|{\bf i}|=N-r.$ This association establishes 
a bijection between such partitions and $I(n,k)$, or, equivalently, 
the Schubert varieties $\Omega_{{\bf i}}$ in ${\Bbb G}_{n,k}$. It is sometimes to convenient 
to denote the Schubert variety $\Omega_{{\bf i}}$ by $\Omega_\nu$ 
where $\nu $ corresponds to ${\bf i}$. This is consistent with our 
notation for a special Schubert variety.

Taking the special Schubert classes 
$[\Omega_i], 1\leq i\leq n-k$, as  algebra 
generators of $H^*({\Bbb G}_{n,k};\bz)$, the structure constants  
are determined by (i) the Pieri formula, 
which expresses the cup-product of an arbitrary 
Schubert class with a special Schubert class as a linear 
combination of with non-negative integral linear combination  
of Schubert classes, and, (ii) the 
Giambelli formula, which expresses an 
arbitrary Schubert class as a determinant in the 
special Schubert classes \cite[Chapter 14]{fulton}.
In particular, the special Schubert classes 
form a set of algebra generators of $H^*({\Bbb G}_{n,k};\bz)$. 
Indeed, $[\Omega_i]=c_i(\gamma^\perp_{n,k})=\bar{c}_i, 
1\leq i\leq n-k$. 

The basis $\{[\Omega_{{\bf i}}]\mid {\bf i}\in I(n,k)\}$ is `self-dual' 
under the Poincar\'e duality. That is, 
assume that ${\bf i},{\bf j}\in I(n,k)$ are such that 
$|{\bf i}|+|{\bf j}|=N$. Then 
$$\langle[\Omega_{{\bf i}}][\Omega_{{\bf j}}],\mu_{n,k}\rangle=\delta_{{\bf  i}',{\bf j}},$$ 
where ${\bf i}'=(n+1-i_k,\cdots, n+1-i_1)\in I(n,k)$. 

The {\it degree} of a Schubert variety 
$\Omega_{{\bf i}}$ of (complex) dimension $r$ is  defined as the 
integer $\langle[\Omega_{{\bf i}}]\bar{c}_1^r,\mu_{n,k}\rangle\in \bz$.  It is well-known \cite{hp},\cite{fulton} that 
$$\deg(\Omega_{{\bf i}})=\frac{r!\prod_{1\leq t<s\leq k}(i_s-i_t)}{(i_1-1)!\cdots (i_k-1)!}\eqno{(1)}$$
In particular 
$$\deg({\Bbb G}_{n,k})=\langle \bar{c}_1^N,\mu_{n,k}\rangle=\frac{N! 1!\cdots (k-1)!}{(n-k)!\cdots (n-1)!}.\eqno{(2)}$$

More generally,  $\deg([\Omega_{{\bf i}}][\Omega_{{\bf j}}])
:=\langle [\Omega_{{\bf i}}][\Omega_{{\bf j}}]\bar{c}_1^q,\mu_{n,k}\rangle
=q!|1/(i_r+j_{k+1-j}-n-1)!|$ where $q=\dim (\Omega_{{\bf i}})+\dim (\Omega_
{{\bf j}})-\dim {\Bbb G}_{n,k}$ (See \cite[p.274]{fulton}.  
We caution the reader that our notations for Grassmann 
manifolds and Schubert varieties are different 
from those used in Fulton's book \cite{fulton}.)

One has the following geometric interpretation for the degree of a Schubert variety.  More generally, given 
any algebraic imbedding  $X\hookrightarrow \bp^m$ of 
a projective variety $X$ of dimension $d$ in the complex 
projective space $\bp^m$, the {\it degree} of $X$ is the    
number of points in the intersection of $X$ with $d$
hyperplanes  in general position. The degree of a Schubert variety 
defined above is the degree of  the Pl\"ucker imbedding 
$\Omega_{{\bf j}}\subset {\Bbb G}_{n,k}\hookrightarrow \bp(\Lambda^k(\bc^n))$, 
defined as  
$U\mapsto \Lambda^k(U)$, where $\Lambda^k(U)$ denotes the $k$-th exterior power of the vector space $U$.  

\noindent
{\bf Cohomology of  Quaternionic Grassmann manifolds}\\
In the case of quaternionic Grassmann manifold ${\Bbb H}{\Bbb G}_{n,k}$, 
one has a Schubert cell decomposition with 
cells only in dimensions $4j$, $0\leq j\leq N$, labeled 
by the same set $I(n,k)$ as in the case of the complex Grassmann manifold $\bc {\Bbb G}_{n,k}$. Furthermore, 
denoting the quaternionic Schubert variety 
corresponding to ${\bf i}\in I(n,k)$ by 
$\Omega_{{\bf i}}^{{\Bbb H}}$,
the structure 
constants defining the integral cohomology algebra of 
${\Bbb H}{\Bbb G}_{n,k}$ for the basis $\{\Omega^{\Bbb H}_{\bf i}\}$ are identical to those in the case of ${\Bbb C}{\Bbb G}_{n,k}$. 
Thus,  the association 
$[\Omega_{{\bf i}}]\mapsto [\Omega_{{\bf i}}^{{\Bbb H}}]$ defines an isomorphism of {\it rings} 
$H^*(\bc{\Bbb G}_{n,k};\bz)\lr H^*({\Bbb H}{\Bbb G}_{n,k};\bz)$ which doubles the degree. In particular 
one has the {\it identical} formula, namely (1),  
for the degrees of quaternionic Schubert  classes.  The orientation on ${\Bbb H}{\Bbb G}_{n,k}$ is chosen 
so the image of the positive 
generator of $H^{2N}(\bc{\Bbb G}_{n,k};\bz)$ under 
the above isomorphism is positive. 

\section{Maps from cohomologically K\"ahler manifolds}

{\it In this section the symbol $d$ will have a 
different meaning from what it did in \S1.}

Let $f\colon X\lr Y$ be any 
continuous map  between two compact connected 
oriented manifolds of the same dimension.
It is well-known that if $f^*$ has non-zero degree, then 
the induced map $f^*\colon H^r(Y;\bz)\lr H^r(X;\bz)$ 
is split-injective for all $r$. In particular, $f^*\colon H^*(Y;\bq)\lr H^*(X;\bq)$  
is a monomorphism of {\it rings}.   

Recall that a compact connected orientable smooth  manifold $X$ is called $c$-{\it symplectic}  (or cohomologically 
symplectic)  if there exists an element $\omega\in H^2(X;\br)$, called a 
$c$-symplectic class,  such that 
$\omega^d\in H^{2d}(X;\br)\cong \br$ is non-zero where 
$d=(1/2)\dim_\br X$. If $\omega$ 
is a $c$-symplectic class in $X$, then $(X,\omega)$ is said to satisfy 
the {\it weak Lefschetz} (respectively {\it hard Lefschetz}) condition if $\cup \omega^{d-1}\colon H^1(X;\br)
\lr H^{2d-1}(X;\br)$ (respectively $\cup \omega^i\colon 
H^{d-i}(X;\br)\lr H^{d+i}(X;\br), 1\leq i\leq d,$) is an isomorphism.
If $(X,\omega)$ satisfies the hard Lefschetz condition, then 
$X$ is called $c$-K\"ahler or cohomologically K\"ahler. 
If $(X, \omega)$ is $c$-K\"ahler, and if $\omega$ is in the image 
of the natural map 
in $H^2(X;\bz)\lr H^2(X;\br)$, we call $X$ $c$-Hodge. 
Note that if $(X,\omega)$ is $c$-K\"ahler and if $H^2(X;\br)\cong \br$, 
then $(X,t\omega)$ is $c$-Hodge for some $t\in \br$.   

Clearly K\"ahler manifolds are  
$c$-K\"ahler and smooth projective varieties over $\bc$ are $c$-Hodge.  
It is known that $\bp^2\#\bp^2$ is $c$-symplectic but not symplectic 
(hence not K\"ahler) 
since it is known that it does not admit even an almost complex 
structure. It is also $c$-K\"ahler. 
Examples of  $c$-symplectic manifolds which satisfy the 
weak Lefschetz condition but not $c$-K\"ahler are also known (cf. \cite{lo}). 

Any $c$-symplectic manifold $(X,\omega)$ is naturally oriented;
the fundamental class of  $X$ will be denoted by $\mu_X\in H_{2d}(X;\bz)\cong \bz$. 

Let $(X,\omega)$ be a $c$-K\"ahler manifold of dimension $2d$. 
Let $1\leq r\leq d$. One has a bilinear form  $(\cdot, \cdot)_\omega$ (or simply $(\cdot,\cdot)$ when there is no danger of confusion) 
on $H^{r}(X;\br)$ defined as $(\alpha, \beta)_\omega=\langle \alpha\beta
\omega^{d-r},\mu_X\rangle, ~\alpha,\beta\in H^r(X;\br)$.  
When $(X,\omega)$ is $c$-Hodge, 
the above form is rational, that is, it restricts to a bilinear form $H^r(X;\bq)\times H^r(X;\bq)\lr \bq$.   It will be important for us to consider 
the bilinear form on the rational vector space $H^r(X;\bq)$ rather than on the real vector space $H^r(X;\br)$.

The bilinear form $(\cdot,\cdot)$ is symmetric (resp. skew symmetric) 
if  $r$ is even (resp. odd). 
Note that the above form is non-degenerate for all $r$. This follows 
from Poincar\'e duality and the hard Lefschetz condition 
that $\beta\mapsto \beta\cup \omega^{d-r}$ is an isomorphism 
$H^r(X;\bq)\lr H^{2d-r}(X;\bq)$.   Further, if $r\leq d$, 
the monomorphism $\cup \omega\colon 
H^{r-2}(X;\bq)\lr H^r(X;\bq)$ is an isometric 
imbedding, i.e., $(\alpha, \beta)=(\alpha\omega,\beta\omega)$ for 
all $\alpha,\beta\in H^{r-2}(X;\br)$.  

As in the case of K\"ahler manifolds 
(cf. \cite{hodge},\cite{weil},\cite{hirz}), one obtains an orthogonal  
decomposition of the real cohomology groups of a $c$-K\"ahler manifold $(X,\omega).$  The decomposition, which preserves the 
rational structure when $(X,\omega)$ is $c$-Hodge,  is obtained 
as follows: Let $1\leq r\leq d$. Let $\mathcal{V}^r_\omega$, or more briefly $\mathcal{V}^r$ when $\omega$ is clear from the context,  be the kernel of the homomorphism 
$\cup \omega^{d-r+1}\colon H^r(X;\br)\lr H^{2d-r+2}(X;\br)$. 
An element of $\mathcal{V}^r$ will be called a {\it primitive class}. 
One has the {\it Lefschetz decomposition}  
$$H^r(X;\br)=\bigoplus_{0\leq q\leq [r/2]}\omega^q \mathcal{V}^{r-2q}.\eqno{(3)}$$ 

We have the following lemma.
\blem   \label{degreehodge}
Suppose that $(X,\omega)$ is a  
$c$-Hodge manifold of dimension $2d$ with second Betti  number equal to $1$.
Let $f\colon X\lr Y$ be any continuous map of non-zero 
degree where $Y$ is a compact manifold with 
non-vanishing second Betti number.  Then:\\
(i) $(\cdot,\cdot)_{t\omega} = t^{d-r}(\cdot,\cdot)_\omega$ on 
$H^{r}(X;\bq)$ for $t\in \bq$, $t\neq 0$.\\
(ii) $(Y,\phi)$ is $c$-Hodge where 
$\phi\in H^2(Y;\bq)$ is the unique class such that $f^*(\phi)=\omega$. Furthermore, $f^*$ preserves the Lefschetz decomposition {\em (3)}, that is, 
$f^*(\mathcal{V}_\phi^r)\subset \mathcal{V}_\omega^r$ for $r\leq d$.\\
(iii) If $\alpha,\beta\in H^r(Y;\bq)$, then $( f^*(\alpha),f^*(\beta))_\omega=\deg(f) (\alpha,\beta)_\phi$. In particular, degree of $f$ equals $\frac{\langle \omega^d,\mu_X\rangle}{\langle \phi^d,\mu_Y\rangle}. $
\elem

\noindent
{\it Proof.}  
(i) This is trivial.\\
(ii) Let $\dim (X)=2d$. 
Since $\deg(f)\neq 0$, $f^*\colon H^{i}(Y;\bq)\lr H^{i}(X;\bq)$ is a 
monomorphism for all $i\leq 2d$. 
Comparing the second Betti numbers of $X$ and $Y$ we conclude 
that $f^*\colon H^2(Y;\bq) \lr H^2(X;\bq)\cong \bq$ is an isomorphism.
Let $\phi\in H^2(Y;\bq)$ be the unique class such that $f^*(\phi)=\omega$.  
Since $f^*$ is a homomorphism of  {\it rings}, we have
$0\neq \omega^d
=(f^*(\phi))^d=f^*(\phi^d)$ and so $\phi^d\neq 0$. 

Let $r\leq d$ be a positive integer. One has a commuting diagram:
$$\begin{array}{ccc}
H^r(Y;\bq)&\stackrel{\cup \phi^{d-r}}{\lr}& H^{2d-r}(Y;\bq)\\
f^*\downarrow ~& &                       ~\downarrow  f^*\\
H^r(X;\bq)&\stackrel{\cup\omega^{d-r}}{\lr}&H^{2d-r}(X;\bq)\\
\end{array}$$

The vertical maps are monomorphisms since $\deg(f)\neq 0$. 
By our hypothesis on $X$, the  homomorphism $\cup \omega^{d-r}$ in the above diagram is 
an isomorphism.  This implies that $\cup\phi^{d-r}$ 
is a monomorphism. Since, by Poincar\'e duality, the vector spaces $H^r(Y;\bq)$ and $H^{2d-r}(Y;\bq)$ have the same dimension, $\cup \phi^{d-r}$
is an {\it isomorphism} and so $(Y;\phi)$ is $c$-Hodge. 
It is clear that $f^*(\mathcal{V}_\phi^r)\subset \mathcal{V}_\omega^r$.  

(iii)  Suppose that $\alpha,\beta\in H^r(Y;\br)$. Then 

$$\begin{array}{lll}(f^*(\alpha), f^*(\beta))_\omega
&=& \langle f^*(\alpha)f^*(\beta)\omega^{d-r};\mu_X\rangle\\
&= &\langle f^*(\alpha\beta)f^*(\phi^{d-r});\mu_X\rangle \\
&=&\langle f^*(\alpha\beta\phi^{d-r});\mu_X\rangle\\
&=&\langle \alpha\beta\phi^{d-r},f_*(\mu_X)\rangle\\
&=&\deg(f)\langle \alpha\beta\phi^{d-r};\mu_Y\rangle\\
&=&\deg(f) (\alpha,\beta)_\phi.\\
\end{array}$$

The formula for the degree of $f$  follows from what has just been established 
by taking $\alpha=\beta =\phi$. \hfill $\Box$

Observe that the summands in the Lefschetz decomposition (3) are mutually 
orthogonal with respect to the bilinear form $(\cdot,\cdot)$.  
Indeed, let $\alpha \in \mathcal{V}^{r-2p},\beta\in \mathcal{V}^{r-2q}, p<q$.
Thus $\alpha\omega^{n-r+2p+1}=0$ and so $\alpha\omega^{n-r+p+q}=0$. 
Therefore  $(\omega^p\alpha,\omega^q\beta)=\langle \alpha\beta\omega^{n-r+p+q},\mu_X\rangle=0.$  As observed earlier the form $(\cdot,\cdot)$ 
is non-degenerate. It follows that the form restricted to each summand 
in (3) is non-degenerate. In favourable situations, the form is 
either positive or negative definite as we shall see in Proposition \ref{mainpropo} below. 

We shall recall some basic results from Hodge theory and 
use several facts concerning harmonic forms, all of which 
can be found in \cite[\S15]{hirz}. They will be needed in the 
proof of Proposition \ref{mainpropo}.

Suppose that $X$ has been endowed with a K\"ahler metric 
with K\"ahler class $\omega\in H^2(X;\br).$
Recall that  one has the decomposition 
$H^r(X;\bc)\cong  \oplus_{p+q=r}
H^{p,q}(X;\bc)$ where $H^{p,q}$ denotes the $\bar{\partial}$-cohomology. 
We identify the $H^{p,q}(X;\bc)$ with the space of harmonic forms 
(with respect to the K\"ahler metric) $B^{p,q}$ 
of type $(p,q)$.  

We shall follow the notations used in \cite[\S 15.8]{hirz}.
One has the operators $L$ and $\Lambda$ on $H^{p,q}(X;\bc)$  
where $L\colon H^{p,q}(X;\bc)\lr H^{p+1,q+1}(X;\bc)$  
equals wedging with the K\"ahler class $\omega$ and $\Lambda \colon 
H^{p,q}(X;\bc)\cong 
B^{p,q}\lr B^{p-1,q-1}\cong H^{p-1,q-1}(X;\bc)$  
is the operator $(-1)^{p+q}\#L\#$ on $B^{p,q}(X;\bc).$ 
The operator $\Lambda$ is dual to $L$ with respect to the 
hermitian scalar product denoted $(\cdot,\cdot)_*$: 
$$(\alpha,\beta)_*:=\int_X\alpha\wedge \#\beta \eqno{(4)}$$ 
on $H^r(X;\bc)=\oplus _{p+q=r}B^{p,q}$.
 
The kernel of $\Lambda$ is denoted by $B^{p,q}_0$.  One has 
the Hodge decomposition 
$$H^{p,q}(X)=\bigoplus_{0\leq k\leq \min\{p,q\}} B_k^{p,q} \eqno{(5)}$$ 
where $B_k^{p,q}:=L^k(B_0^{p-k,q-k})$ 
is the space of all harmonic forms $\phi$ of 
type $(p,q)$ and class $k$. Then  
the distinct summands in (5) are pairwise orthogonal with respect 
to $(\cdot,\cdot)_*$.
Also, $\Lambda L^k$ is a non-zero scalar multiple of $L^{k-1}$  on 
$B^{p-k,q-k}_0$ for $p+q\leq d, 1\leq k\leq \min\{p,q\}$.   

\bpropo\label{mainpropo} 
Suppose that $(X,\omega)$ is a compact connected  
K\"ahler manifold such that $H^{p,q}(X;\bc)
=0$ for $p\neq q$. Then the form $(-1)^{q+r}(\cdot,\cdot)_\omega$ restricted to $\omega^q\mathcal{V}^{2r-2q}\subset H^{2r}(X;\br)$ is positive definite for 
$0\leq q\leq r, 1\leq r\leq [d/2]$. 
\epropo

\noindent
{\it Proof.} 
First assume that $d=\dim_\bc X$ is even, say $d=2s$.
In view of  our hypothesis, all odd Betti numbers of $X$ vanish and we have  
$B_k^{p,q}=0$ for all $p\neq q, k\geq 0,$ so that 
$$H^{2r}(X;\bc)=H^{r,r}(X;\bc)=\oplus_{0\leq k\leq r}B_k^{r,r}.\eqno{(6)}$$

The real cohomology group $H^{2r}(X;\br)\subset 
H^{2r}(X;\bc)=H^{r,r}(X;\bc)$ has an orthogonal decomposition induced from (3):
$$H^{2r}(X;\br)=\oplus_{0\leq k\leq s}E_k^{r,r}\eqno{(7)}$$ 
where $E_k^{p,p}=\{\alpha\in B_k^{p,p}\mid \alpha=\bar{\alpha}\}$. 
Now taking $r=s=d/2$ one has  
$\#\alpha=(-1)^{s+k}\alpha$ for $\alpha\in E_k^{s,s}$. In 
particular the bilinear form (4) equals $(-1)^{s+k}Q$ where 
$Q(\alpha,\beta)=\int_X \alpha \beta.$   Therefore $(-1)^{s+k}Q$ restricted to each $E_k^{s,s}$ is positive definite. 
 
We shall show in Lemma \ref{primitiveclasses} below that 
$\omega^k\mathcal{V}^{d-2k}=E_k^{s,s}$.   
The proposition follows immediately from this since  
$(\alpha,\beta)=(\omega^{s-r}\alpha,\omega^{s-r}\beta)$ 
for $\alpha,\beta\in \omega^k\mathcal{V}^{2r-2k}$ as $d=2s,$ completing the proof in this case.   

Now suppose that $d$ is odd. Consider the 
K\"ahler manifold $Y=X\times \bp^1$ where  
we put the Fubini-Study 
 metric on $\bp^1$ with K\"ahler class $\eta$ 
being the `positive' generator of $H^2(\bp^1;\bz)\subset H^2(\bp^1;\br) $  
and the product structure on $Y$ so that the 
K\"ahler class of $Y$ equals $\omega+\eta=:\phi$. 
By K\"unneth theorem $H^*(Y;\br)=
H^*(X;\br)\otimes H^*(\bp^1;\br)$.  We shall identify the cohomology 
groups of $X$ and $\bp^1$ with their images in $H^*(Y;\br)$ via the 
monomorphisms induced by the first and second projection respectively. 
Under these identifications, 
$H^{p,q}(Y;\bc)= H^{p,q}(X;\bc)\oplus H^{p-1,q-1}(X;\bc)\otimes H^{1,1}(\bp^1;\bc)$.   In particular, $H^{p,q}(Y;\bc)=0$ unless $p=q$. By what has been proven already, the form $(-1)^{r+k}(\cdot,\cdot)$ is positive definite on 
$\phi^k\mathcal{V}^{2r-2k}_\phi\subset H^{2r}(Y;\br)$.   

Choose a base point in $\bp^1$ and consider the inclusion map $j
\colon X\hookrightarrow Y$.  
The imbedding $j$ is dual to $\eta$. Also $j^*(\phi)=
\omega$. It follows that $j^*(\phi^k\mathcal{V}^{2r-2k}_\phi)\subset 
\omega^k\mathcal{V}^{2r-2k}_\omega$ 
for $0\leq k<r,~1\leq r<d$.  Since the kernel of $j^*\colon H^{2r}(Y;\br)\lr H^{2r}(X;\br)$ equals $H^{2r-2}(X;\br)\otimes H^2(\bp^1;\br)$, and maps 
$H^{2r}(X;\br)\subset H^{2r}(Y;\br)$ isomorphically onto $H^{2r}(X;\br)$, 
we must have $j^*(\phi^k\mathcal{V}^{2r-2k}_\phi)=\omega^k\mathcal{V}^{2r-2k}_\omega$.

Let $\alpha,\beta\in H^{2r}(X;\br)\subset H^{2r}(Y;\br)$.  
Since $j\colon X\hookrightarrow Y$ is dual to $\eta$, we have $j_*(\mu_X)=\eta\cap\mu_Y$. Therefore,  
 $$\begin{array}{lll}
(j^*(\alpha),j^*(\beta))_\omega&=&\langle j^*(\alpha\beta)j^*(\omega)^{d-2r};\mu_X \rangle\\
        &=&\langle \alpha\beta \omega^{d-2r},j_*(\mu_X)\rangle\\
        &=&\langle \alpha\beta\omega^{d-2r},\eta\cap \mu_Y\rangle \\
        &=&\langle \alpha\beta\omega^{d-2r}\eta,\mu_Y\rangle\\
 \end{array}$$
Since $\eta^2=0$ we have $\phi^{d-2r+1}=\omega^{d-2r+1}+(d-2r+1)\omega^{d-2r}\eta$.  Furthermore, $\alpha\beta\omega^{d-2r+1}\in H^{2d+2}(X;\br)
=0$.  Therefore, we conclude that 
$(j^*(\alpha),j^*(\beta))_\omega=\frac{1}{d-2r+1}\langle \alpha\beta\phi^{d-2r+1},\mu_Y\rangle
=\frac{1}{d-2r+1}(\alpha,\beta)_\phi$.
This shows that the bilinear form $(\cdot,\cdot)_\omega$ on $H^{2r}(X;\br)$ is a positive 
multiple of the form $(\cdot,\cdot)_\phi$ on $H^{2r}(Y;\br)$ restricted to $H^{2r}(X;\br)$. 
It follows that the bilinear form $(-1)^{r+k}(\cdot,\cdot)$ on $H^{2r}(X;\br)$ restricted to $\omega^k\mathcal{V}^{2r-2k}(X)$ is positive definite.\hfill $\Box$        
 
We must now establish the following 
  
\blem \label{primitiveclasses} With notations as above, 
assume that $d=2s$ is even. Under the hypothesis 
of the above proposition, 
$E_k^{s-k,s-k}$ equals $\omega^k\mathcal{V}^{d-2k}$, $0\leq k\leq s$.  
\elem
{\it Proof.} Since $L$ preserves real forms, it suffices to show that  
$E_0^{r,r}=\mathcal{V}^{2r}$ when $r\leq s$.  By definition 
$E_0^{r,r}=B_0^{r,r}\cap H^{2r}(X;\br)=
\{\alpha\in H^{r,r}(X;\bc)\mid \Lambda(\alpha)=0,~\alpha=\bar{\alpha}\}.$  

Let $\alpha\in E_0^{r,r}$. Suppose that $p\geq 1$ 
is the largest integer such that 
$\omega^{d-2r+p}\alpha=:\theta$ 
is a {\it non-zero} real harmonic form of type $(d-r+p,d-r+p)$. Since 
$L^{d-2r+2p}\colon H^{r-p,r-p}(X;\bc)\lr H^{d-r+p,d-r+p}(X;\bc)$  is an isomorphism, and since $\omega$ is 
real   
there must be a real form $\beta\in H^{r-p,r-p}(X;\br)$ such that $L^{d-2r+2p}(\beta)=\theta
=L^{d-2r+p}(\alpha)$. 
Since $p$ is the largest, using the decomposition (6) we see that $\beta\in B_0^{r-p,r-p}$. 
 Applying $\Lambda^{d-2r+p}$ both sides and 
(repeatedly) using $\Lambda L^q\beta$ is a {\it non-zero} 
multiple of $L^{q-1}\beta $ when $r-p+q<d$ we see that 
$\beta$ is a non-zero multiple of  $\Lambda^{p}\alpha=0$.  
Thus $\beta=0$ and 
hence $\theta=0$, which contradicts our assumption.  Therefore 
$L^{d-2r+1}(\alpha)=0$ and so $\alpha\in \mathcal{V}_0^r$.  On the other 
hand $\Lambda $ maps $H^{2r}(X;\bc)$ onto $H^{2r-2}(X;\bc)$. 
A dimension argument shows that $E_0^{r,r}=
\mathcal{V}^{2r}$. \hfill $\Box$                        


\beg\label{grassvanish}
The Grassmann manifold ${\Bbb G}_{n,k}$ has the structure of 
a K\"ahler manifold with K\"ahler class $\omega:=\bar{c}_1=[\Omega_1]\in H^2({\Bbb G}_{n,k};\bz)$.  (This fact follows, for example, from the Pl\"ucker imbedding 
${\Bbb G}_{n,k}\hookrightarrow \bp^{{{n}\choose{k}}-1}$. )
 The bilinear form $(\cdot,\cdot)$ is understood to be defined with 
 respect to $\omega$.
 An orthogonal basis for $\mathcal{V}^{2r}_{n,k}\subset H^{2r}({\Bbb G}_{n,k};\bq)$ can be obtained inductively using Gram-Schmidt
 orthogonalization process as follows. Recall from \S 2 the basis  $\bar{\mathcal{C}}_r$ 
 for $H^{2r}({\Bbb G}_{n,k};\bq)$. Clearly $\omega\cdot\bar{\mathcal{C}}_{r-1}=\bar{c}_1\cdot \bar{\mathcal{C}}_{r-1}=
 \{\bar{c}^{{\bf j}}\in \bar{\mathcal{C}}_r\mid j _1>0\}$ is a basis for $\omega H^{2r-2}({\Bbb G}_{n,k};\bq)$. Therefore we see that the subspace spanned by 
 $\bar{\mathcal{C}}_{r,0}:=\{\bar{c}^{{\bf j}}
 \in \bar{\mathcal{C}}_{r}\mid j_1=0\}$ is  complementary to 
 $\oplus_{q>0}B_q^{r-q,r-q}\subset H^{2r}
({\Bbb G}_{n,k};\bq).$  The required basis is obtained by taking the orthogonal 
projection of $\bar{\mathcal{C}}_{r,0}$ onto $\mathcal{V}^{2r}$. 
Indeed, inductively assume that an orthogonal 
basis $\{v_{{\bf j}}\}$ for $\omega H^{2r-2} ({\Bbb G}_{n,k};\bq)$ that is compatible with the direct sum decomposition $\oplus _{q>0}B_q^{r-q,r-q}$
has been constructed. 
We need only apply the orthogonalization process to the (ordered) set 
$\{v_{{\bf j}}\}\cup \{\bar{c}^{\bf j}\in \bar{\mathcal{C}}_r\mid j_1=0\}$ where the 
elements $\bar{c}^{{\bf  j}}$ are ordered, say, according to 
lexicographic order 
of the exponents.  For example, taking $n=12, k=6, r=6,$ the elements of 
$\bar{\mathcal{C}}_{6,0}$  are ordered as $\bar{c}_2^3,\bar{c}_2\bar{c}_4,\bar{c}_3^2,\bar{c}_6.$ We denote 
the basis element of $\mathcal{V}^{2r}$ obtained from $c^{{\bf j}}\in \bar{\mathcal{C}}_{r}$ by  $v_{{\bf j}}$. Note that when $r\leq k$,  the span of the set  
$\{ v_{{\bf j}}\mid j_r=0\}\subset H^{2r}({\Bbb G}_{n,k};\bq)$  equals space of the 
decomposable elements in $H^{2r}({\Bbb G}_{n,k};\bq)$ since, according to 
our assumption on the ordering of elements $\bar{c}^{{\bf j}}$, the 
element $\bar{c}_r$ is the greatest and so $v_r$ does not occur 
in any other $v_{{\bf j}}$.  Thus $v_r-\bar{c}_r$ belongs to the 
ideal  $\mathcal{D}\subset H^*({\Bbb G}_{n,k};\bq)$ and $v_{{\bf j}}\in \mathcal{D}$ 
for all other ${\bf j}.$   

We illustrate this for $r=2,3$. (When $r=1$, 
 $\mathcal{V}^1=0$. )
The element $v_2=\bar{c}_2-\frac{(\bar{c}_2,\omega^2)}{(\omega,\omega)}\omega^2=
\bar{c}_2-\frac{\deg \bar{c}_2}{\deg {\Bbb G}_{n,k}}\omega^2\in 
H^4({\Bbb G}_{n,k};\bq)$ 
is a basis for the one-dimensional space $\mathcal{V}^4$.

Similarly, $v_3$ is a basis for $\mathcal{V}^6$ where 
$$\begin{array}{lll}
v_3&:=&\bar{c}_3-
\frac{(\bar{c}_3,v_2\omega)}{(v_2\omega,v_2\omega)}v_2\omega -\frac{(\bar{c}_3,\omega^3)}{(\omega^3,\omega^3)}\omega^3\\
&=&\bar{c}_3-\frac{\deg \bar{c}_3}{\deg {\Bbb G}_{n,k}}\omega^3-
\frac{\deg {\Bbb G}_{n,k}\deg(\bar{c}_3\bar{c}_2 )-\deg \bar{c}_2\deg\bar{c}_3}
{\deg {\Bbb G}_{n,k}\deg(\bar{c}_2^2)-(\deg\bar{c}_2)^2}v_2\omega.\\
\end{array}$$ 

This leads to \\
$(v_3,v_3)=(v_3,\bar{c}_3)=
\deg (\bar{c}_3^2)-\frac{(\deg \bar{c}_3)^2}{\deg {\Bbb G}_{n,k}}
-\frac{\deg(\bar{c}_3\bar{c}_2 )\deg {\Bbb G}_{n,k}-\deg \bar{c}_2\deg\bar{c}_3}
{\deg {\Bbb G}_{n,k}\deg(\bar{c}_2^2)-(\deg\bar{c}_2)^2}\deg(\bar{c}_3v_2).$
\eeg
 
The following calculation will be used in the course of the 
proof of Theorem \ref{degreezero}.
 
\blem \label{v2v2}
With the above notation, 
$(v_2,v_2)=  \deg {\Bbb G}_{n,k}\frac{(k^2-1)((n-k)^2-1)}{2(N-1)^2(N-2)(N-3)}.$

\elem 
\noindent
{\it Proof.}  The proof involves straightforward but lengthy calculation which we work out below.

Since $(v_2,\bar{c}_1^2)=0$,  we get 
$(v_2,v_2)=(v_2,c_2)
=(\bar{c}_2,\bar{c}_2)-\frac{\deg\bar{c}_2}{\deg {\Bbb G}_{n,k}}(\bar{c}_2,\omega^2)
=\deg {\Bbb G}_{n,k}(\frac{\deg(\bar{c}_2^2)}{\deg {\Bbb G}_{n,k}}-(\frac{\deg\bar{c}_2}{\deg {\Bbb G}_{n,k}})^2).$

Since  
$\bar{c}_2^2=[\Omega_2]^2=[\Omega_4]+[\Omega_{3,1}]+[\Omega_{2,2}]$, we see that   
$\frac{\deg \bar{c}_2^2}{\deg {\Bbb G}_{n,k}}= \frac{\deg\bar{c}_4}{\deg {\Bbb G}_{n,k}}+
\frac{\deg\Omega_{3,1}}{\deg {\Bbb G}_{n,k}}+\frac{\deg\Omega_{2,2}}{\deg 
{\Bbb G}_{n,k}}.\\$

Now an explicit calculation yields, upon using $N=k(n-k)$:
$$\begin{array}{lll}
\frac{\deg \bar{c}_4}{\deg {\Bbb G}_{n,k}}&=&\frac{(n-k-1)(n-k-2)(n-k-3)(k+1)(k+2)(k+3)}{4!(N-1)(N-2)(N-3)},\\
\frac{\deg\Omega_{3,1}}{\deg {\Bbb G}_{n,k}}
&=&\frac{(n-k+1)(n-k-1)(n-k-2)(k+2)(k+1)(k-1)}{2!4(N-1)(N-2)(N-3)},\\
\frac{\deg\Omega_{2,2}}{\deg {\Bbb G}_{n,k}}
&=&\frac{N(k-1)(k+1)(n-k+1)(n-k-1)}{2! 3\cdot 2(N-1)(N-2)(N-3)},\\
\frac{\deg \bar{c}_2}{\deg {\Bbb G}_{n,k}}
&=&\frac{(k+1)(n-k-1)}{2!(N-1)}.
\end{array}$$

Substituting these in the above expression for $(v_2,v_2)$ 
we get 
$(v_2,v_2)=\frac{(k+1)(n-k-1)}{4!(N-1)^2(N-2)(N-3)}A$ 
where, again using $N=k(n-k)$ repeatedly, \\ 
$$\begin{array}{ll}
A&:=(N-1)\{(n-k-2)(k+2)(n-k-3)(k+3)\\
&+3(n-k-2)(k+2)(n-k+1)(k-1)
+2N(k-1)(n-k+1)\}\\
&-6(N-2)(N-3)((n-k-1)(k+1))^2)\\
&=(N-1)\{(N+2(n-2k)-4)(N+3(n-2k)-9)\\
&+3(N+2(n-2k)-4)(N-(n-2k)-1)+2(N-(n-2k)-1)\}\\
&-6(N-2)(N-3)(N+(n-2k)-1)\\
&=12(N-(n-2k)-1)\\
&=12(k-1)(n-k+1).\\
\end{array}$$
Therefore, 
$(v_2,v_2)=\deg {\Bbb G}_{n,k}\frac{(k^2-1)((n-k)^2-1)}{2(N-1)^2(N-2)(N-3)}.$\hfill $\Box$

\brem\label{quaternion} 
Although quaternionic Grassmann manifolds 
are not $c$-K\"ahler,  one could use the symplectic 
Pontrjagin class $\eta:=e_1(\gamma_{n,k})\in H^4({\Bbb H}
{\Bbb G}_{n,k};\bz)$ in the place of $\bar{c}_1\in H^2
({\Bbb C}{\Bbb G}_{n,k};\bz)$ to define a 
pairing $(\cdot,\cdot)_\eta$ on $H^{4r}({\Bbb H}{\Bbb 
G}_{n,k};\bq)$ and the primitive classes $v_j\in 
H^{4j}({\Bbb H}{\Bbb G}_{n,k};\bq)$. 
We define $\mathcal{V}^{4r}\subset H^{4r}({\Bbb H}{\Bbb G}_{n,k};\bq)$ to be the kernel of 
$$\cup\eta^{N-2r+1}
\colon H^{4r}({\Bbb H}{\Bbb G}_{n,k};\bq)
\lr H^{4N-4r+4}({\Bbb H}{\Bbb G}_{n,k};\bq).$$ 
The form $(\cdot,\cdot)_\eta$ is definite when restricted 
to the space $\eta^q\mathcal{V}^{4r-4q}\subset H^{4r}({\Bbb H}{\Bbb G}_{n,k};\bq)$. 
The formula given in Lemma \ref{v2v2} holds 
without any change.  These statements follow 
from the degree doubling isomorphism from the 
cohomology algebra of ${\Bbb G}_{n,k}$ to that 
of ${\Bbb H}{\Bbb G}_{n,k}$ which maps 
the $i$-th Chern class of the tautological complex 
$k$-plane bundle over ${\Bbb G}_{n,k}$ to the $i$-th symplectic Pontrjagin 
class of the tautological left ${\Bbb H}$-bundle 
over ${\Bbb H}{\Bbb G}_{n,k}$. 
\erem

\section{Proofs of  Main Results}
In this section we prove the main results of the paper, namely 
Theorems \ref{degdet}, \ref{degreezero} and \ref{degreeone}.  {\it We will only consider the 
case of complex Grassmann manifolds.} The 
proofs in the case of quaternionic Grassmann manifolds 
follow in view of the fact that the cohomology 
algebra of ${\Bbb H}{\Bbb G}_{n,k}$ is 
isomorphic to that of ${\Bbb C}{\Bbb G}_{n,k}$ 
via an isomorphism that doubles the degree. 

Recall that 
complex Grassmann manifolds are smooth projective varieties and that 
Schubert subvarieties yield an {\it algebraic} cell decomposition. In  
particular their Chow ring is isomorphic to singular cohomology (with $\bz$-coefficients) via an isomorphism that doubles the 
degree. It follows that $H^{p,q}({\Bbb G}_{n,k};\bc)=0$ 
for $p\neq q$.  Therefore results of the previous section hold for ${\Bbb G}_{n,k}$. The bilinear form $(\cdot,\cdot)$ is understood to be 
defined with respect to $\omega=\bar{c}_1\in H^2({\Bbb G}_{n,k};\bz)\cong \bz$. 

\blem\label{degree}
Let $f\colon {\Bbb G}_{n,k}\lr {\Bbb G}_{m,l}$ be any continuous 
map where $k(n-k)=l(m-l)$.  Suppose that $f^*(c_1(\gamma^\perp_{m,l}))=\lambda c_1(\gamma^\perp_{n,k})$ where $\lambda\in \bz$. Then 
$$\deg(f) =\lambda^N \frac{\deg {\Bbb G}_{n,k}}{\deg{\Bbb G}_{m,l}}.$$   
\elem

\noindent 
{\it Proof.} This follows immediately from 
Lemma \ref{degreehodge}(i) and (iii). \hfill $\Box$

\noindent
{\it Proof of Theorem \ref{degreeone}.} 
We may suppose that ${\Bbb F}=\bc$ and that $l\leq k$; otherwise  
$k\leq l\leq [m/2]$ in which case $\deg(f)=0$ for any $f$  by 
\cite[Theorem 2]{rs}.   

Suppose that $\deg(f)=\pm 1$. We have  
$$\begin{array}{rcl}
\frac{\deg {\Bbb G}_{n,k}}{\deg {\Bbb G}_{m,l}}
&=&\frac{1!\cdots (k-1)!(m-l)!\cdots (m-1)!}{1!\cdots (l-1)! (n-k)!\cdots (n-1)!}\\
~&=&\frac{l! \cdots (k-1)!(m-l)!\cdots (m-1)!}{(n-k)!\cdots (n-1)!}\\
~&=&(\prod _{1\leq j\leq k-l}\frac{(l-1+j)!}{(n-k+j-1)!})(\prod_{1\leq j\leq l}\frac{(m-j)!}{(n-j)!})\\
\end{array}$$

Note that after simplifying 
$(l+j-1)!/(n-k+j-1)!$ for each $j$ in the first product, we 
are left with product of $(k-l)$ blocks of 
$(n-k-l)$ consecutive positive integers in the 
denominator, the {\it largest} to occur being $(n-l-1)$. Similar 
simplification in the second product yields a 
product of $l$ blocks of $(m-n)$ consecutive 
integer, the {\it smallest} to occur being $(n-l+1)$. 
Since $(k-l)(n-k-l)=l(m-n)$ we conclude that 
$\deg({\Bbb G}_{n,k})>\deg({\Bbb G}_{m,l})$. 

In the notation of Lemma \ref{degree} above,  
we see that either $\deg(f)=0$ or $|\deg(f)|>|\lambda|^N\geq 1$---
a contradiction. Therefore $(m,l)=(n,k)$ if $\deg(f)=\pm 1$. 
Now $f^*\colon H^*({\Bbb G}_{n,k};\bz)\lr H^*({\Bbb G}_{n,k};\bz)$ induces an isomorphism.  
Since ${\Bbb G}_{n,k}$ is a 
simply connected CW complex, by Whitehead's theorem, $f$ is a homotopy equivalence. \hfill $\Box$ 

\brem
(i) The above is a topological analogue of the 
result of Paranjape and Srinivas \cite{parsri} that any non-constant 
morphism\\ $f\colon {\Bbb G}_{n,k}\lr {\Bbb G}_{m,l}$ is an {\it isomorphism} of 
varieties provided the ${\Bbb G}_{m,l}$ is not the projective space.  
Our conclusion in the topological realm is weaker. 
Indeed it is known that there exist continuous self-maps of any complex and quaternionic 
Grassmann manifold which have large positive degrees.  See \cite{friedlander} and also \cite{shiga}.  

(ii) Endomorphisms of the cohomology algebra of 
${\Bbb G}_{n,k}$ have been classified by M. Hoffman \cite{hoffman}. 
They are either `grading homomorphisms' defined by 
$c_i\mapsto \lambda^i c_i, 1\leq i\leq k$ for some $\lambda$ 
or when $n=2k$, the composition of a grading homomorphism 
with the homomorphism induced by the diffeomorphism  
$\perp\colon {\Bbb G}_{n,k}\lr {\Bbb G}_{n,k}$ defined as $U\mapsto U^\perp$. 
\erem

Recall from Example \ref{grassvanish} the construction of the primitive classes $v_j\in H^{2j}({\Bbb G}_{n,k};\bq), 2\leq j\leq k$. To avoid possible confusion,  
we shall denote the primitive classes in $H^{2j}({\Bbb G}_{m,l};\bq)$ corresponding to $j=2,\cdots, l$ by $u_j$. Also $\mathcal{V}_{m,l}^{2r}\subset H^{2r}({\Bbb G}_{m,l};\bq)$ will denote the space of primitive classes. 
The following lemma is crucial for the proof of Theorem \ref{degdet}.

\blem\label{crucial} Suppose that $f\colon {\Bbb G}_{n,k}\lr {\Bbb G}_{m,l}$ is 
a continuous map such that  $f^*(c_1(\gamma_{m,l}^\perp))=\lambda c_1(\gamma_{n,k}^\perp)=\lambda\bar{c}_1$ 
with $\lambda\neq 0$. Let $2\leq j\leq l$. Then, with the above notations,  
$f^*(u_j)=\lambda_j v_j$    
where $\lambda_j\in \bq$ is such that $$\lambda_j^2 =
\lambda^{2j}\frac{\deg {\Bbb G}_{n,k}}{\deg {\Bbb G}_{m,l}}\frac{(u_j,u_j)}{(v_j,v_j)}$$  
for $2\leq j\leq l$. 
\elem
 
\noindent
{\it Proof.}  The degree of $f$ equals $\lambda^N\deg {\Bbb G}_{n,k}/\deg {\Bbb G}_{m,l}\neq 0$ by Lemma \ref{degree}.

Therefore $f^*\colon H^{2j}({\Bbb G}_{m,l};\bq)\lr H^{2j}({\Bbb G}_{n,k};\bq )$ 
is an isomorphism and $f^*(\mathcal{V}^{2j}_{m,l})=\mathcal{V}_{n,k}^{2j}$, 
since $f^*$ is a monomorphism and the dimensions are equal as $j\leq l$.    
Note that $f^*$ maps the space of decomposable elements $\mathcal{D}_{m,l}^{2j}\subset H^{2j}({\Bbb G}_{m,l};\bq)$  
isomorphically onto $\mathcal{D}_{n,k}^{2j}$.  
Since $u_j\perp \mathcal{D}_{m,l}^{2j}\cap \mathcal{V}_{m,l}^{2j}$ we 
see that, by Lemma \ref{degreehodge} (ii), $f^*(u_j)\perp \mathcal{D}_{n,k}^{2j}\cap \mathcal{V}_{n,k}^{2j}$. As the form $(\cdot, \cdot)$ on $\mathcal{V}^{2j}_{n,k}$ is definite by Proposition 
and $\mathcal{V}^{2j}_{n,k}=\bq v_{j}\oplus (\mathcal{V}^{2j}_{n,k}\cap \mathcal{D}^{2j}_{n,k})$ is an orthogonal decomposition, we must have 
$f^*(u_j)=\lambda_j v_j$ for some $\lambda_j\in \bq$. 

Recall that 
$\deg(f)=\lambda^N\deg {\Bbb G}_{n,k}/\deg {\Bbb G}_{m,l}$.  Note that 
$$\begin{array}{lll}
\lambda^{N-2j}(f^*(u_j),f^*(u_j))&=&
(f^*(u_j),f^*(u_j))_{\lambda\bar{c}_1}\\ 
&=&\deg(f) (u_j,u_j)_\omega\\ 
& = & \lambda^N \frac{\deg {\Bbb G}_{n,k}}{\deg {\Bbb G}_{m,l}} (u_j,u_j)\\
\end{array}$$ 
by Lemma \ref{degreehodge}. 
Thus $\lambda_j^2(v_j,v_j)=(f^*(u_j),f^*(u_j))=\lambda^{2j}\frac{\deg {\Bbb G}_{n,k}}{\deg {\Bbb G}_{m,l}}(u_j,u_j)$.  \hfill $\Box$

We are now ready to prove Theorem \ref{degdet}.

\noindent
{\it Proof of Theorem \ref{degdet}}: We need only 
consider the case ${\Bbb F}=\bc$.
  Recall that the cohomology 
algebra $H^*({\Bbb G}_{m,l};\bz)$ is generated by $\bar{c}_1,
,\cdots, \bar{c}_l$ where $\bar{c}_j=c_j(\gamma_{m,l}^\perp)$. 
Therefore $f^*\colon H^*({\Bbb G}_{m,l};\bz) \lr H^*({\Bbb G}_{n,k};\bz)$ is determined by the images of $\bar{c}_j, 
1\leq j\leq l$. 

As observed in Example \ref{grassvanish}, 
one has $u_j-\bar{c}_j\in \mathcal{D}_{m,l}^{2j}, 2\leq j\leq l$.  
It follows easily by induction that each $\bar{c}_j, 1\leq j\leq l,$ 
can be expressed as a polynomial with 
rational coefficients in $\bar{c}_1, u_2,\cdots,u_l$.
Therefore $\bar{c}_1=:u_1,u_2, \cdots, u_l$ generate $H^*({\Bbb G}_{m,l};\bq)$. 

Lemma \ref{degree} implies that $f^*(u_1)=\lambda c_1(\gamma_{n,k}^\perp)$ where $\lambda^N$---and hence $\lambda$ upto a sign---is determined by the degree of $f$.

Now by Lemma \ref{crucial}, the image of $f^*(u_j)=\lambda_jv_j$ 
where $\lambda_j$ is determined upto a sign. \hfill $\Box$ 

\noindent
{\it Proof of Theorem \ref{degreezero}:}  We assume,  
as we may that, ${\Bbb F}=\bc$. 
We  preserve the 
notations used in the above proof. 
Recall from Lemma \ref{v2v2} that 
$(v_2,v_2)=\deg {\Bbb G}_{n,k} \frac{(k^2-1)((n-k)^2-1)}{2
(N-1)^2(N-2)(N-3)}. $ 
Therefore, by Lemma \ref{crucial}
we have 
$$\begin{array}{lll}
\lambda^2_2&=&\lambda^4\frac{\deg {\Bbb G}_{n,k}}{\deg {\Bbb G}_{m,l}}\frac{(v_2,v_2)}{(u_2,u_2)}\\
&=&\lambda^4 (\frac{\deg {\Bbb G}_{n,k}}{\deg {\Bbb G}_{m,l}})^2
\frac{(k^2-1)((n-k)^2-1)}{(l^2-1)((m-l)^2-1)}\\
&=& B^2(k^2-1)(l^2-1)((n-k)^2-1)((m-l)^2-1)
\end{array}$$
where $B:=\frac{\lambda^2\deg {\Bbb G}_{n,k}}{\deg {\Bbb G}_{m,l}(l^2-1)
((m-l)^2-1)}\in \bq$. 
It follows that $\deg(f)=0$ unless $Q:=(l^2-1) (k^2-1)((m-l)^2-1)((n-k)^2-1)$ is a perfect square.
It 
remains to show that there are at most finitely many values for $m,n$ 
for which the $Q$ is a perfect square. This is proved in following proposition.

\bpropo \label{perfectsquare}
Let $1<a<b$ be positive integers. Then 
there are at most finitely many solutions in $\bz$ for the system of equations 
$$y^2=Q(a,b,x,z),  ~az=bx,\eqno{(8)}$$ where  
$Q(a,b,x,z):=(a^2-1)(b^2-1)(x^2-1)(z^2-1)$. 
\epropo
\noindent
{\it Proof.}  Let $r=\textrm{gcd}(a,b)$ and write 
$a=rs, b=rt$ so that $tx=sz$.  Then the system of 
equations (8) can be rewritten as 
$ y^2=F(x)$ where $F(x):=(1/s^2)(a^2-1)(b^2-1)(x^2-1)(t^2x^2-s^2).$ 
Note that $F(x)\in \bq[x]$ has distinct zeros in $\bq$. By a theorem of 
Siegel  \cite[Theorem D.8.3, p. 349]{hs} it follows that 
the equation $y^2=F(x)$ has only finitely many solutions 
in the ring $R_S\subset K$ of $S$-integers where $K$ 
is any number field and $S$ any finite set of absolute 
valuations of $K$, including all archimedean valuations. In particular, taking $K=\bq$ and $S$ the usual (archimedean) 
absolute value, we see that there are only finitely many integral solutions of (8). \hfill $\Box$
   
For the rest of the paper we shall only be concerned 
with the number theoretic question of $Q(a,b,c,d)$ 
being a perfect square. 

\brem
(i) We observe that there are {\it infinitely} many 
integers $1<a<b<c<d$ such that $Q(a,b,c,d)$ is a perfect square. Indeed given $a,b,$ let 
$c$ be any positive integer such that   
$(a^2-1)(b^2-1)(c^2-1)=Pu^2$ where $P>1$ is square 
free. Let $(x,y)$ be any 
solution with $x\neq 0$ of the so called Pell's equation $y^2=1+Px^2$. Then $d=|y|$ is a solution whenever $d>c$. Since the Pell's equation has infinitely many solutions, 
there are infinitely many such $d$.  

(ii) Suppose that $(l^2-1)(k^2-1)(c^2-1)=x^2$ is 
a perfect square.  (There exists such  
positive integers $c$---in fact infinitely many of them--- 
for which this happens if and only if $(l^2-1)(k^2-1)$ is 
{\it not} a perfect square.) Then there does {\it not} exist any $d>1$ such that $Q(l,k,c,d)$ is a perfect square. 
Assume further that $l|(kc)$---this can be arranged, for 
example, taking $k$ to be a multiple of $l$---and set $n:=c+k$, $m:=kc/l$ so 
that $k(n-k)=l(m-l)$. Then $Q(l,k,n-k,m-l)$ is not 
a perfect square.

(iii) We illustrate below situations  $Q(l,k,n-k,m-l)$ 
is not a perfect square (assuming that $k(n-k)=l(m-l)$) depending on congruence classes modulo a suitable prime power of the parameters involved. \\
(1) For an odd prime $p,$ suppose that $k\equiv 
p^{2r-1}\pm 1\mod p^{2r}$ and none of the numbers $l,m-l,n-k$ is congruent to $\pm 1 \mod p$. Then  
$p^{2r-1}|Q$ but $p^{2r}\not| Q$. \\
(2) Suppose that  $m\equiv l\equiv 5 \mod 8,$ and $k\equiv 7\mod 16$.   Then $(m-l)^2-1$ is odd, $l^2-1
\equiv 8 \mod 16$, $k^2-1\equiv 16\mod 32$ and 
$l(m-l)=k(n-k)$ implies $(n-k)$ is even and so $(n-k)^2-1$ is  odd. Thus $Q\equiv 2^7\mod 2^8$.\\
(3) Suppose that $l\equiv 0\mod 8, m\equiv l \mod 2$, 
$k\equiv 3 \mod 8$ , then $Q\equiv 8\mod 16$.
\erem

We conclude the paper with the following  
 
\bpropo   \label{examples}
Let $c>1$ and let $k=3$ or $7$.  Suppose that  
$Q(2,k,2c,kc)$ is a perfect square. Then there exists integers $\xi,\eta, v>1$ such that   
$c=\frac{1}{2}(\xi^2\eta^2+1), ~\xi^2\eta^2-3v^2=-2$
and (i) $\xi^2-3\eta^2=-2$ when $k=3$ and (ii) $\xi^2-7\eta^2=-6$ 
when $k=7$. 
\epropo
\noindent
{\it Proof.} Assume that $k=7$ and that    
$Q:=Q(2,7,2c,7c)=3^2 2^4(2c-1)(2c+1)(7c-1)(7c+1)$ is a perfect square. 
There are several cases to consider depending on the gcd of the pairs of numbers involved. 
Write $(2c-1)=\alpha u^2, 2c+1=\beta v^2, 7c-1=\gamma x^2, 7c+1=\delta y^2,$ where $\alpha, \beta,\gamma, \delta$ are square free integers. 
Since $Q$ is a perfect square and since 
$\gcd(2c-1,2c+1)=1,\gcd(7c-1,7c+1)=1 ~\textrm{or} ~2, \gcd(2c\pm 1,7c\pm 1)=1,~ \textrm{or}~ 5, \gcd(2c\pm 1,7c\mp 1)=1,3,~\textrm{or} ~9$, 
the possible values for  $(\alpha,\beta)$ are:
$(1,1),  (1,5),  (1,3), (3,1), (5,1),
(1,15),(15,1),(5,3),(3,5)$. 
The possible values for $(\gamma,\delta)$ are the same as 
for $(\alpha,\beta)$ as well as $(2\alpha,2\beta)$. 

Suppose $(\alpha,\beta)=(1,1)$. Since $(2c-1)+2=(2c+1)$, 
we obtain $u^2+ 2=v^2$ which has no solution. 
If $(\alpha,\beta)=(3,1), $ then $3u^2+2=v^2$. This equation 
has no solution mod $3$. Similar arguments show that if $(\alpha,\beta)=(5,1),(1,5),(1,15), (15,1), (5,3)$, there are no solutions for $u,v$.  If $(\alpha,\beta)=(3,5)$, then $(\gamma,\delta)=(5,3)$ or $(10,6)$. If 
$(\gamma,\delta)=(5,3)$ again there is no solution  mod $3$ for the equation $5x^2+2=3y^2$. 
When $(\gamma,\delta)=(10,6)$ we obtain $10x^2+2=6y^2$. 
This has no solution mod $5$.  

It remains to consider the case $(\alpha,\beta)
=( 1,3)$. In this case obtain the equation $u^2+2=3v^2$ which has 
solutions, for example, $(u,v)=(5,3)$. 
Now $(\alpha,\beta)=(1,3)$ implies $(\gamma,\delta)=(3,1)$ or $(6,2)$. 
If $(\gamma,\delta)=(3,1)$ then we obtain the equation $3x^2+2=y^2$
which has no solution mod $3$.  
So assume that $(\gamma,\delta)=(6,2)$. 
As $(\alpha,\delta)=(1,2)$ we obtain 
$4y^2-7u^2=9$, that is,  
$4y^2-7u^2=9$. Thus $(2y-3)(2y+3)=7u^2$. Either 
$7|(2y-3)$ or $7|(2y+3)$. Say $7|(2y-3)$ and write $(2y-3)=7z$. Now 
$z(7z+6)=u^2$.  Observe that $\gcd(z,7z+6)$ divides $6$. 
Since $\beta=3$, $2c-1=u^2$ is not divisible by $3$. Also, $u$ being odd, we must have $\gcd(z,7z+6)=1$.  It follows that both $z, 7z+6$  
are perfect squares. This forces that $6$ is a square mod $7$---a contradiction.  
Finally, suppose that $7|(2y+3)$. Then repeating 
the above argument we see that both $(2y-3)=:\eta^2$ and $(2y+3)/7=:\xi^2$ are 
perfect squares. It follows that $7\xi^2-6=\eta^2$ is  a perfect square. Hence $2c-1=u^2=\xi^2\eta^2$. Since $2c+1=3v^2$, the proposition follows.

We now consider the case $k=3$. We merely sketch the proof in this 
case.  Let, if possible, $Q$ = $2^33(2c-1)(2c+1)(3c-1)(3c+1)$ be a perfect square. 
Write  $2c-1 = \alpha u^2, 2c+1 = \beta v^2, 3c-1 = \gamma x^2, 3c+1 = \delta y^2$, where $\alpha,\beta,\gamma,\delta$ are square free integers and $u,v,x,y$ are positive integers.
Arguing as in the case $k=7$, 
following are the only possible values for $\alpha, \beta, \gamma,\delta$: 
$(\alpha , \beta )=(1,3),(3,1),(3,5),(5,3),(1,15),(15,1)$, and $(\gamma , \delta )=(1,2),(2,1),(2,5),(5,2),(1,10),(10,1)$. 
It can be seen that only the case 
$( \alpha , \beta , \gamma , \delta ) = (1,3,2,1)$ remains to be 
considered, the remaining possibilities leading to contradictions. 
Thus we have  $2c-1= u^2 , 2c+1 =3v^2, 3c-1 =2 x^2$ and $3c +1 = y^2$. Therefore, we have $4 x^2 -1  = 3 u^2$, i.e., $(2x-1)(2x+1) = 3 u^2$ . Hence, $3|(2x-1)$ or $3|(2x+1)$. 

Suppose that $3|(2x-1)$.  Write $3z = 2x-1$, $z \in {\Bbb Z} $. 
Since $z$ is odd,  we have 
$\gcd(z, 3z+2)= 1$. As $z(3z+2)=u^2$ we conclude that 
$z$ and $3z+2$ have to be perfect squares. 
This implies that $2$ is a quadratic residue mod $3$--a contradiction. Therefore $3\not|(2x-1)$ and we must have $3|(2x+1)$ and 
both $z$ and $3z-2$ will have to be perfect squares. Write $z = \eta ^2 $ and $3z-2 = \xi ^2 $ so that $\xi ^2 - 3 \eta ^2 = -2 $ and  $v^2=u^2+2=\xi^2\eta ^2+2.$ This completes the proof.\hfill $\Box$

\brem
(i) Let $K=\bq[\sqrt{7}]$ and let $R$ be the ring of integers 
in $K$. If $\xi+\eta\sqrt{7}\in R$, then $\xi,\eta\in \bz$. 
Denote the multiplicative ring of units in $R$ by $U$. 
Note that any element of $U$ has norm $1$. (This is 
because $-1$ is a quadratic non-residue mod $7$.) 
Using Dirichlet Unit theorem $U$ has rank $1$; indeed  
$U$ is generated by $\nu:=(8+3\sqrt{7})$ and $\pm 1$. 
The integers $\xi,\eta$ as in the above proposition yield an element  
$\xi+\eta\sqrt{7}$ of norm $-6$ and the set $S\subset R$ of all elements of norm $-6$ is stable under the multiplication action by $U$.  An easy argument 
shows that $S$ is the union of orbits through 
$\lambda :=1+\sqrt{7},\bar{\lambda}=1-\sqrt{7}$. 
Thus $S=\{\pm\lambda \nu^k,\pm\bar{\lambda}\nu^k\mid k\in\bz\}$. 
  
Observe that if $ \xi,\eta$ are as in Proposition \ref{examples}(ii), 
then $\xi+\sqrt{7}\eta\in S$.  Listing elements $\xi+\eta\sqrt{7}\in  S$ with $ \xi,\eta>1$ in increasing order of 
$\eta$, the first three elements are 
$13+5\sqrt{7}, 29+11\sqrt{7}, 209+79\sqrt{7}$. 
Straightforward verification shows that when $\xi+\eta\sqrt{7}$ is  
equals any of these, then there does not exist an integer $v$ such that 
$\xi^2\eta^2+2=3v^2$.  Since the next term is $463+175\sqrt{7}$, we have the lower bound  
$2c>175^2 \times 463^2=6565050625$ in order that $Q(2,7,2c,7c)$ be a perfect square (assuming $c>1$).  

(ii) Now, let $K=\bq[\sqrt{3}]$  and let $R$ be the ring of integers in $K$. Note that $\xi+\eta\sqrt{3}\in R$, then $\xi,\eta\in{\Bbb Z}$. Denote the multiplicative ring of units in $R$ by $U$, which is generated by 
generated by $(2+\sqrt{3})$ and $\pm 1$. 

Suppose that $Q(2,3,2c,3c)$ is a perfect square, $c>1$. Then 
the integers $\xi,\eta$, as in the above proposition, yield an element $\xi +\eta \sqrt{3}$ of norm $-2$. 
The set $S \subset R$ of all elements of norm $-2$ is stable under the multiplication action by $U$. In fact it can be verified easily that 
$S=\{\pm(1+\sqrt{3})(2+\sqrt{3})^m\mid m\in\bz\}$. 

Listing these with $\xi , \eta > 1$, in increasing order of $\eta$, the first five elements  are  $5+3\sqrt{3},19+11\sqrt{3}, 71+41\sqrt{3},265+153\sqrt{3},  989+571\sqrt{3}.$ If $\xi+\eta\sqrt{3}$ 
equals any of these, direct verification shows that there is no integer $v$ satisfying the equation $\xi^2\eta^2+2=3v^2$.
The next term of the sequence being $3691 + 2131 \sqrt{3}$  we obtain the lower bound $2c > 2131^2 \times 3691^2
=61866420601441$.   
\erem

Perhaps, the above arguments can be applied for 
possibly other values of $k$ and $l$ to obtain 
lower bounds for $n$, particularly when $k+l$ and $k-l$ are primes or prime powers. 

\noindent
{\bf Acknowledgments:}  We are grateful to 
R. Balasubramanian and D. S. Nagaraj for  
their help with Proposition \ref{perfectsquare}. We 
thank Nagaraj also for his valuable comments.

\end{document}